\documentclass[e-only,10pt,reqno]{ofj}

\usepackage{setspace}
\usepackage{lineno}
\usepackage{color}
\usepackage[dvipsnames,svgnames,x11names]{xcolor}
\usepackage[hang,flushmargin]{footmisc}
\usepackage{orcidlink}
\usepackage{cite}

\usepackage[
nonumberlist, 
acronym,      
toc,          
]{glossaries}
\newacronym{amr}{AMR}{adaptive mesh refinement}
\newacronym{bc}{BC}{boundary condition}
\newacronym{bdf}{BDF}{backward differentiation formula}
\newacronym{bosss}{BoSSS}{bounded support spectral solver}
\newacronym{cfl}{CFL}{Courant-Friedrichs-Lewy}
\newacronym{dg}{DG}{discontinuous Galerkin}
\newacronym{dggerman}{DG}{Diskontinuierliche Galerkin}
\newacronym{dof}{DOF}{degrees of freedom}
\newacronym{eoc}{EOC}{experimental order of convergence}
\newacronym{fem}{FEM}{finite element method}
\newacronym{fe}{FE}{finite element}
\newacronym{fvm}{FVM}{finite volume method}
\newacronym{fv}{FV}{finite volume}
\newacronym{fds}{FDS}{finite differencing scheme}
\newacronym{ic}{IC}{initial condition}
\newacronym{ip}{IP}{interior penalty}
\newacronym{ivp}{IVP}{initial value problem}
\newacronym{bvp}{BVP}{boundary value problem}
\newacronym{ode}{ODE}{ordinary differential equation}
\newacronym{pde}{PDE}{partial differential equation}
\newacronym{rhs}{RHS}{right-hand side}
\newacronym{roc}{ROC}{rate of convergence}
\newacronym{sip}{SIP}{symmetric interior penalty}
\newacronym{rans}{RANS}{Reynolds-averaged Navier--Stokes}
\newacronym{les}{LES}{large-eddy simulation}
\newacronym{dns}{DNS}{direct numerical simulation}
\newacronym{evm}{EVM}{eddy-viscosity model}
\newacronym{nlevm}{NLEVM}{nonlinear eddy-viscosity model}
\newacronym{rsm}{RSM}{Reynolds stress model}
\newacronym{earsm}{EARSM}{explicit algebraic Reynolds stress model}
\newacronym{rdt}{RDT}{rapid distortion theory}
\newacronym{pdf}{PDF}{probability-density function}
\newacronym{cfd}{CFD}{computational fluid dynamics}
\newacronym{mms}{MMS}{method of manufactured solutions}
\newacronym{cas}{CAS}{computer algebra system}
\newacronym{blas}{BLAS}{basic linear algebra subprograms}
\newacronym{lapack}{LAPACK}{linear algebra package}
\newacronym{mpi}{MPI}{message-passing interface}
\newacronym{lmn}{LMN}{Lundgren--Monin--Novikov}
\newacronym{sst}{SST}{shear-stress transport}
\newacronym{simple}{SIMPLE}{semi-implicit method for pressure-linked equation}
\newacronym{dae}{DAE}{differential-algebraic equation}

\newacronym{hdf}{HDF}{hierarchical data format}

\newacronym{asbl}{ASBL}{asymptotic suction boundary layer}

\setacronymstyle{short-long}

\usepackage{hyperref}
\usepackage{tikz}
\usetikzlibrary{shadows, shapes.arrows, arrows, shapes, decorations.pathreplacing,3d}

\tikzset{
  ashadow/.style={opacity=.25, shadow xshift=0.07, shadow yshift=-0.07},
}

\tikzstyle{decision_templ} = [diamond, draw, fill=blockbg,
    text width=4.5em, text badly centered, node distance=3cm, inner sep=0pt,
    drop shadow={ashadow}]
\tikzstyle{block_templ} = [rectangle, draw, fill=blockbg,
text width=21em, text centered, minimum height=0em,
    drop shadow={ashadow}]
\tikzstyle{block_templ} = [rectangle, draw, fill=blockbg,
    text width=18em, text centered, minimum height=0em,
    drop shadow={ashadow}]
\tikzstyle{line} = [draw, thick, -latex']
\tikzstyle{cloud_templ} = [draw, ellipse,fill=blockbg, node distance=3cm,
minimum height=2em, text width=8em, text centered,
    drop shadow={ashadow}]
\colorlet{blockbg}{blue!5}

\usetikzlibrary{arrows,shapes,decorations.pathreplacing,3d}
\usetikzlibrary{backgrounds}
\usepackage{latex_macros/dakling}
\toggletrue{EinsteinNotationDiff}

\catcode`_=12 %
\renewcommand{\texttt}[1]{%
  \begingroup
  \ttfamily
  \begingroup\lccode`~=`/\lowercase{\endgroup\def~}{/\discretionary{}{}{}}%
  \begingroup\lccode`~=`[\lowercase{\endgroup\def~}{[\discretionary{}{}{}}%
  \begingroup\lccode`~=`.\lowercase{\endgroup\def~}{.\discretionary{}{}{}}%
  \begingroup\lccode`~=`_\lowercase{\endgroup\def~}{_\discretionary{}{}{}}%
  \catcode`/=\active\catcode`[=\active\catcode`.=\active\catcode`_=\active
  \scantokens{#1\noexpand}%
  \endgroup
}
\catcode`_=8 %

\let\paragraph\undefined

\copyrightinfo{2022}{The Authors. This work is an open access article under the {\color{blue}{\href{https://creativecommons.org/licenses/by-sa/4.0/}{CC BY-SA 4.0}}} license}
\newcommand{\OF}[0]{OpenFOAM\textsuperscript{\textregistered}}
\newcommand{\of}[0]{\OF}
\newcommand{\foamdg}[0]{\code{foam-dg}}
\newcommand{\mono}[0]{\code{mono}}
\newcommand{\microsoft}{Microsoft\textsuperscript{\textregistered}}

\newcommand{\textcite}[1]{\cite{#1}}
\newcommand{\parencite}[1]{\cite{#1}}

\OpenFOAMversions{foam-extend 4.1}

\ifdefined\review
\Repository{https://github.com/FDYdarmstadt/foam-dg}
\else
\Repository{https://github.com/FDYdarmstadt/foam-dg}
\fi

\DOI{TBD} 

\ifdefined\review
\fi

\begin{document}


\title[Coupling \of{} with BoSSS]{Coupling \of{} with BoSSS, a discontinuous Galerkin
  solver written in \csharp{}}

\ifdefined\review
  \author{}
  \address{}
  \email{}
\else



\author{Dario Klingenberg$^{1}$\orcidlink{0000-1234-5678-0000}}
\author{Hrvoje Jasak$^{3}$}
\author{Holger Marschall$^{2,4}$}
\author{Florian Kummer$^{1,2,*}$}
\address{$^1$Chair of Fluid Dynamics, TU Darmstadt, Darmstadt, Germany}
\email{kummer@fdy.tu-darmstadt.de}
\address{$^2$Graduate School Computational Engineering, TU Darmstadt, Darmstadt, Germany}
\address{$^3$Department of Physics, University of Cambridge, Cambridge, UK}
\address{$^4$Mathematical Modeling and Analysis, TU Darmstadt, Darmstadt, Germany}
\fi


\date{\today}

\dedicatory{}

\maketitle

\ifdefined\review
  \linenumbers
\else
\fi

\begin{abstract}
  In this article, we present the \foamdg{} project,
  which provides a bridge
  between \of{} and the high-order \gls{dg} framework \gls{bosss}.
  Thanks to the flexibility of the coupling approach, mixed calculations where
  some parts of the equation system are solved in \of{} and others are solved in
  \gls{bosss} are easily possible.
  This is showcased using the convective Cahn--Hilliard equation, where the
  Cahn--Hilliard part is solved in \gls{bosss} and the Navier--Stokes part is
  solved in \of{}.
  The obtained results appear reasonable, though the main focus of this paper is
  to present and document the \foamdg{} project rather than on quantitative
  results.
\end{abstract}

\section{Introduction}
\label{sec:introduction}

Recent years have seen a growing interest within the \gls{cfd} community in
high-order methods.
While the \gls{fvm} implemented in \of{} can realize high convergence orders by
increasing the cell stencil on structured grids, practical applications with
complex geometries often rely on automatically generated unstructured grids.
Moreover, in parallel calculations, the increased stencil size results in
a significant communication overhead.
The \gls{dg} method offers the attractive possibility to accomplish high-order
convergence even on unstructured grids, and it can generally handle relatively
low-quality grids well.
Unlike in \gls{fvm}, the \gls{dg} method approximates the solution in a cell
with a polynomial of order \(p\). It can be shown that the resulting L\(^2\)-norm
convergence order is \(p+1\), allowing for arbitrarily high convergence order
without increasing the cell stencil size. A short introduction into the method
is given in \cref{sec:glsdg-method}.

Several \gls{dg} solvers are currently available.  Popular \cpp{}-frameworks
include Nektar++ \parencite{nektarpp2015} and deal.II \parencite{dealII94} which
both feature continuous and discontinuous Galerkin \gls{fem} solvers.
The highly modular \cpp{}-library Dune \parencite{dunefemdg:17}, among many things,
also features a \gls{dg} solver.
A spectral element method solver with a strong focus on performance is provided
by NEK5000 \parencite{nek5000}, which is written in \clang{} and \fortran{}.
In \textcite{xu2017}, the authors describe implementing a \gls{dg} method on top
of \of{}. However, their resulting code HopeFOAM is still in an early stage of
development.

In this article, we present a method of enhancing \of{} with the high-order
\gls{dg} framework \gls{bosss}.
This combination of codes bases poses a number of challenges, not least because
\of{} is written in the programming language \cpp{}, whereas \gls{bosss} is
implemented in \csharp{}. This is discussed in detail in \cref{sec:comm-betw-c++}.
Nevertheless, \gls{bosss} offers a number of unique features.  The interactive
BoSSSpad allows not only for interactive pre- and postprocessing, but also
permits the prototyping and implementation of entire solvers through a Jupyter
notebook interface. Even though this is not immediately relevant to
the presently discussed \of{}-interface because calculations are started from
\of{}, we expect this to significantly flatten the learning curve for users that
later decide to pursue a full high-order implementation in \gls{bosss}.
Furthermore, \gls{bosss} is written in a modern object-oriented style and
features a continuous integration system using a state-of-the-art testing
framework.

Apart from the computer science aspects arising from the combination of
programming language runtimes, there is also the question of how to translate
between the \gls{fvm} and the \gls{dg} world in a sensible way. This question is
addressed in \cref{sec:coupling-fvm-dg}.

As shown in figure \labelcref{fig:flow},
the example considered here is the solution of the Cahn--Hilliard equation with
convection. The momentum equation is solved in \of{}, the resulting fields
are then communicated to \gls{bosss}, where the Cahn--Hilliard equation is
solved. After synchronizing the resulting fields back to \of{}, this loop
continues from the beginning.
The rationale behind this particular way of coupling, apart from serving as an
interesting setup showcasing the modularity of our approach, is that
high-order accuracy is much more beneficial for the Cahn--Hilliard equation, whereas
for the solution of the momentum equation, a standard \gls{fvm} may
suffice in this context.
The reason why the Cahn--Hilliard equation benefits more from being solved using
a high-order method lies in its unusual characteristics, such as a cubic term
and a fourth derivative.

Despite the runtime performance penalty arising as a consequence of this code
coupling, the approach discussed here offers a number of interesting advantages.
In particular, existing low-order \of{} solvers can be extended to use
a high-order discretization with relatively little effort while also relying on
a large existing \gls{bosss} functionality and features.
This makes this approach particularly attractive for prototyping the extension
of existing \of{} solvers with high-order methods.

\newcommand*{\cdg}[1]{C^{#1}_\text{DG}}
\newcommand*{\cfvm}[1]{C^{#1}_\text{FVM}}
\newcommand*{\udg}[1]{U^{#1}_\text{DG}}
\newcommand*{\ufvm}[1]{U^{#1}_\text{FVM}}

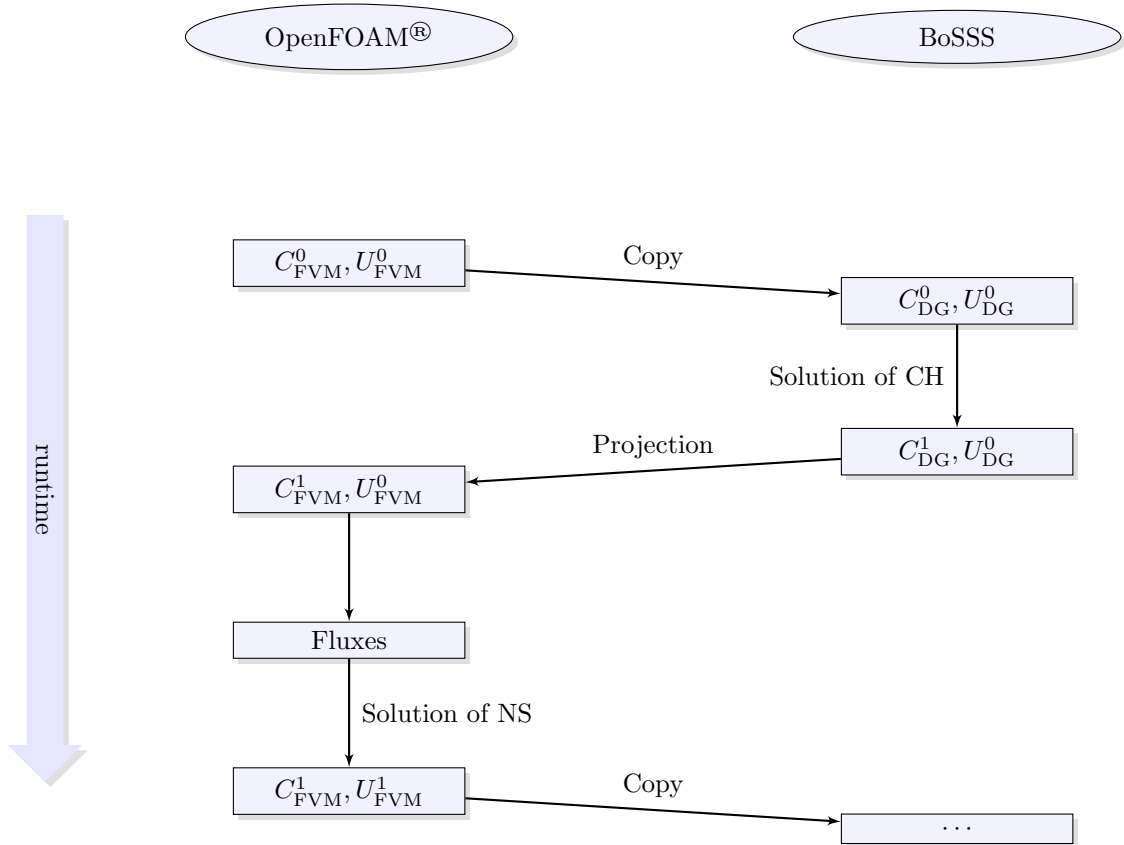
\begin{figure}
\begin{tikzpicture}
\tikzstyle{block} = [block_templ, text width=8em]
\tikzstyle{smallblock} = [block_templ, text width=5em]
\tikzstyle{cloud} = [cloud_templ, text width=8em]
\tikzstyle{smallcloud} = [cloud_templ, text width=5em]

  \node[cloud] (of) at (-4,0) {\of{}};
  \node[cloud] (bosss) at (4,0) {BoSSS};
  \node [fill=blue!10, text width=20em, text centered, single arrow, drop shadow={ashadow}, rotate=-90] at (-8, -6) {runtime};

  \node[block] (c0) at (-4,-3) {\(\cfvm{0}, \ufvm{0}\)};
  \node[block] (c0bosss) at (4,-3.5) {\(\cdg{0}, \udg{0}\)};
  \path[draw, thick, -latex'] (c0) -- node[above=0.1em] {Copy} (c0bosss);


  \node[block] (c1bosss) at (4,-5.5) {\(\cdg{1}, \udg{0}\)};
  \path[draw, thick, -latex'] (c0bosss) -- node[left=0.1em] {Solution of CH} (c1bosss);

  \node[block] (c1of) at (-4,-6.0) {\(\cfvm{1}, \ufvm{0}\)};
  \path[draw, thick, -latex'] (c1bosss) -- node[above=0.1em] {Projection} (c1of);

  \node[block] (c1offlux) at (-4,-8.0) {Fluxes};
  \path[draw, thick, -latex'] (c1of) -- node[above=0.1em] {} (c1offlux);

  \node[block] (c1of1) at (-4,-10.0) {\(\cfvm{1}, \ufvm{1}\)};
  \path[draw, thick, -latex'] (c1offlux) -- node[right=0.1em] {Solution of NS} (c1of1);

  \node[block] (c1bosss) at (4,-10.5) {\(\cdots\)};
  \path[draw, thick, -latex'] (c1of1) -- node[above=0.1em] {Copy} (c1bosss);
\end{tikzpicture}
\label{fig:flow}
\caption{Flow chart of the coupling between \of{} and BoSSS.}
\end{figure}

\subsection{The \gls{dg} method}
\label{sec:glsdg-method}

In the following, we give a brief introduction into the \gls{dg} method aimed at
readers already familiar with \gls{fvm}.
As a simple example, consider the one-dimensional \gls{pde}
\begin{equation}
  \label{eq:7}
  \pd{\phi}{t} + \pd{f(\phi)}{x} = 0,
\end{equation}
where \(\phi\) is a scalar variable and \(f(\phi)\) is some function.
In the first step, \cref{eq:7} is multiplied by a test function \(\test\) and
integrated over the entire domain volume \(\Omega\), leading to the weak form
\begin{equation}
  \label{eq:8}
  \int_{\Omega}{\pd{\phi}{t} \test}dV + \int_{\Omega} \pd{f(\phi)}{x} \test dV = 0.
\end{equation}
Integration by parts and Stokes's theorem allows rewriting this into
\begin{equation}
\label{eq:9}
\int_{\Omega}{\pd{\phi}{t} \test}dV
+ \int_{\Omega}  f(\phi) \pd{\test}{x} dV
- \int_{\partial \Omega} {f(\phi)} \test \vec{n} dS
= 0,
\end{equation}
where \(\partial \Omega\) denotes the boundary of the domain and
\(\vec{n}\) is the normal vector pointing outside of the domain.
We now discretize the domain into cells, where the \(\ith\) cell occupies the
volume \(\Omega_{i}\) and has the boundary \(\partial \Omega_{i}\),
\begin{equation}
  \label{eq:10}
\sum_{i} \int_{\Omega_{i}}{\pd{\phi}{t} \test}dV
+ \sum_{i} \int_{\Omega_{i}}  f(\phi) \pd{\test}{x} dV
- \sum_{i} \int_{\partial \Omega_{i}} {f(\phi)} \test \vec{n} dS
= 0.
\end{equation}

Note that at this stage, the classical \gls{fvm} derivation can still be
recovered by setting \(\test = 1\) and assuming \(\phi\) to be constant across
each cell.
However, in \gls{dg}, we now proceed by making a polynomial ansatz for the trial
function \(\phi\) and the test function \(\test\).
For \(\phi\), cell-local polynomials are used. In the \(\ith\) cell, the ansatz reads
\begin{equation}
  \label{eq:224}
  \phi_{i}(\xv, t) = \sum_{k={0}}^{N_{i}} \hat{\phi}_{i; k}(t) \psi_{i;k}(\xv),
\end{equation}
where \(\psi_{k;i}(\xv)\) is the polynomial basis (e.g.\ monomials, i.e.\
\(\psi_{i;0} = 1, \enskip \psi_{i;1} = x, \enskip \psi_{i;2} = x^{2}, \ldots \),
though other choices, such as Lagrange or Legendre polynomials, are usually more
suitable for practical applications), and \(\hat{\phi}_{i; k}\) are the degrees
of freedom in the cell.
As is the characteristic feature of a Galerkin method, the same basis
functions are used for the test function \(\test\) as for the trial function
\(\phi\), i.e.\
\begin{align}
  \label{eq:dgansatz}
v_{i; k} (\xv) = \psi_{i;k}(\xv).
\end{align}

\newcommand{\phiin}[1]{\phi_{#1}^{\text{in}}}
\newcommand{\phiout}[1]{\phi_{#1}^{\text{out}}}
\newcommand{\phiini}{\phiin{i}}
\newcommand{\phiouti}{\phiout{i}}
\newcommand{\vin}{\test^{\text{in}}}
\newcommand{\vout}{\test^{\text{out}}}
Note that
this implies that in contrast to \gls{fem}, continuity across cells is not
ensured by the ansatz. This poses a problem for the surface integral in
\cref{eq:10}, because \(\phi\) is not uniquely defined on interior cell
boundaries. As a result, this term has to be replaced by a so-called flux
function \(\hat f(\phiini, \phiouti, \vec{n}, \test)\), leading to
\begin{equation}
  \label{eq:11}
\sum_{i} \int_{\Omega_{i}}{\pd{\phi}{t} \test}dV
+ \sum_{i} \int_{\Omega_{i}}  f(\phi) \pd{\test}{x} dV
- \sum_{i} \int_{\partial \Omega_{i}} {\hat f(\phiini, \phiouti, \vec{n}, \test)} dS
= 0.
\end{equation}
A suitable choice for \(\hat f\) reconciles the value inside the cell \(\phiini\)
with that outside the cell \(\phiouti\) (which could be the \(\phiin{i+1}\) of the
neighboring cell or also a boundary condition) while also ensuring that
discontinuities are penalized.
As a result, continuity across cells is weakly enforced, whereas in \gls{fem},
it is strongly enforced.
The most important advantage over classical \gls{fem} in the context of
\gls{cfd} is that an appropriate choice of the flux function can ensure that
conservation laws are locally fulfilled.


\subsection{The Cahn--Hilliard equation}
\label{sec:cahn-hill-equat-1}

First published in \textcite{cahn1958}, the Cahn--Hilliard equation reads
\begin{equation}
  \label{eq:1}
  \pd{c}{t} + \ui{k} \pd{c}{x_{k}} = D \pd[2]{\mu}{x_{k}},
\end{equation}
where the chemical potential \(\mu\) is given by
\begin{equation}
  \label{eq:2}
  \mu = c^{3} - c - \gamma \pd[2]{c}{x_{k}},
\end{equation}
\(U\) is the fluid velocity, \(c\) is the order parameter, which ranges between
negative and positive unity and indicates phase domains, and \(D\) and
\(\gamma\) are empirical parameters.  This equation describes the phenomenon of
phase separation, where one phase is identified with \(c = -1\), and the other
one with \(c = 1\), with intermediate values denoting a mixture of the two phases.

Apart from its physical significance, the Cahn--Hilliard equation is also
interesting from a numerics point of view due to its somewhat unusual features,
which include the cubic nonlinearity in \cref{eq:2} and the fact that it is a
fourth-order differential equation if \cref{eq:1,eq:2} are combined into a
single differential equation.

This part of the equation system is solved in \gls{bosss}. The second-derivative
terms are discretized using the \gls{sip} method
\parencite{arnold1982,arnold2002}.

In the present solver, the velocity field \(\ui{i}\) is obtained from the
solution of the Navier--Stokes equations, which read
\begin{align}
  \label{eq:konti}
    \pd{U_{i}}{x_{i}} &= 0,\\
    \label{eq:nast}
  \pd{\ui{i}}{t}
    + \ui{j} \pd{\ui{i}}{x_{j}}
  &=
   - \pd{P}{x_{i}}
    + \pdl{\left(\nu \pd{\ui{i}}{x_{j}}\right)}{x_{j}},
\end{align}
where \(P\) is the pressure (divided by the density) and \(\nu\) stands for the kinematic viscosity.
As mentioned above, this part of the system is solved in \of{}.

\section{Implementation}
\label{sec:implementation}

In this section, we present the details of the solver implementation. Note that
this implementation involves three software projects (\of{}, \gls{bosss} and
\foamdg{}), as well as two programming languages (\cpp{} and \csharp{}).
Apart from the resulting computer science related challenges, a brief discussion
of numerical aspects is also given.

\subsection{Communication between \cpp{} and \csharp{}}
\label{sec:comm-betw-c++}

The main developments described here take place in the \foamdg{} project,
which acts as a bridge between \of{} and \gls{bosss}. Additionally, some
implementation work is done on the \gls{bosss} side, and minor changes are
made to \of{}. Note that we use the \of{} version foam-extend4.1 in
this work.

While \gls{bosss} usually runs on the \dotnet{} framework developed by \microsoft{}, for
the purpose of coupling with \of{}, the alternative \mono{} framework is used.
Originally developed for running \csharp{} on non-Windows operating
systems, active development has slowed down with the release of the cross-platform
\dotnetfive{} framework. However, \mono{}'s ability to interface with \clang{} (and,
thus, \cpp{}) makes it well-suited for the present purposes.

The embedding approach uses the library \code{libmono} and is documented at
\url{https://www.mono-project.com/docs/advanced/embedding/}.
On a basic level, it offers the possibility for a \clang{}/\cpp{} program to instatiate \csharp{}
classes and call \csharp{} methods. 
For the currently discussed program, this means that a simulation is run like a
normal \of{} calculation, and the code then calls \gls{bosss} functionality internally.
The first step to accomplish this is to compile and link the \mono{} runtime when
compiling \foamdg{}.
Thanks to the \of{} build system, the necessary commands do not have to be
entered by the user each time, but are integrated into the project source files.
The second step is, at runtime, to initialise the \mono{} runtime using the
\code{mono\_jit\_init} function. Subsequently, the \gls{bosss} assembly is loaded
using \code{mono_domain_assembly_open}. The code implementing this is located in
the \file{foam-dg/src/discontinuousGalerkin/BoSSSwrapper/Globals.cpp} file.  On
the \gls{bosss} side, classes and methods to be exposed to \cpp{} have to be
marked with the \code{CodeGenExport} attribute.
This in turn allows the necessary wrapper code on the \cpp{} side to be
automatically generated by \gls{bosss}.

As a result, the development procedure for adding \gls{bosss} functionality to
\of{} is as follows:
(i) implement the desired features in \gls{bosss}, make sure to mark the
relevant classes and methods with the  \code{CodeGenExport} attribute,
(ii) compile \gls{bosss},
(iii) create the necessary bindings on the \foamdg{} side by invoking the
\file{foam-dg/src/discontinuousGalerkin/BoSSSwrapper/createBindings.sh}
shell script, and
(iv) compile \foamdg{} using the \file{foam-dg/src/Allwmake} shell script
and the \file{foam-dg/src/discontinuousGalerkin/BoSSSwrapper/compile.sh} shell
script. This step requires foam-extend4.1 to be installed correctly.
After these steps, the solver is ready to run. For the Cahn--Hilliard
solver, example cases are available in \file{foam-dg/run/}, e.g.\
\code{dropletInShearFlow}.
These cases can be run like a normal \of{} case. For convenience, an
\file{Allrun} shell script performing all the necessary steps is also included
in each example case. After the necessary preparations involving grid creation
using \file{blockMeshDict} and initializing the fields using
\code{funkySetFields} provided by the \code{swak4foam} library, the case is run
by invoking the \code{CahnHilliardFoam} command.
This starts a normal \of{} solver at first, though \gls{bosss} is initialized
in the early part of the run together with other typical set-up operations such
as creating time and the fields.
Note that the \foamdg{} library implements \gls{dg} counterparts for all relevant
\of{} concepts, such as the mesh, the fields and the matrices.
These objects serve to store the data communicated to and from \gls{bosss}.
The \gls{bosss} code enters the picture when the \code{dgScalarMatrix} is
assembled. In the case of this solver, it only contains one term, the
\code{dgm::dgCahnHilliard} operator.
If configured properly in the \file{system/dgSchemes} file, this operator
calls the \code{bosssCahnHilliard} code, which in turn invokes \gls{bosss}.
After creating the exported \gls{bosss} classes \code{OpenFoamMatrix} (which
contains the grid and the fields) and \code{OpenFoamPatchField} (containing the
boundary conditions) on the \of{} side, it calls the static \gls{bosss}
method \code{FixedOperators.CahnHilliard()}.
It is this method that performs the actual core part of solving the
Cahn--Hilliard equation, which is done in \gls{bosss}. After the method
finishes, the results are synchronized to \of{}.
At this point, enough information is available on the \of{} side to solve the
momentum equation, after which point the time step is finished.

Note that it is also possible for testing purposes to run calculations entirely
within \gls{bosss}, though in that case, there is no straightforward method
available for coupling to the momentum equation.
This is valuable during code development and in order to integrate the
\gls{bosss} side of the solver into \gls{bosss}'s continuous testing framework.

\subsection{Implementation details on the \of{} side}
\label{sec:impl-deta-of}

The project \foamdg{} does not implement a \gls{dg} full solver, but it
defines some classes to accommodate the information given to and returned from
\gls{bosss}.

From an \of{} developer's point of view, adding the \gls{bosss} operator is as
simple as using it to define the \code{dgScalarMatrix} (i.e.\ \foamdg{}'s
equivalent to the \code{fvmScalarMatrix}),
\begin{lstlisting}[emph={ddt,div,dgCahnHilliard}]
dgScalarMatrix CEqn(
  dgm::dgCahnHilliard(C, U, Phi, Flux)
);
\end{lstlisting}
and to specify
\begin{lstlisting}
dgCahnHilliardSchemes
{
  dgCahnHilliard(CU) bosssCahnHilliard;
}
\end{lstlisting}
in \file{system/dgSchemes}.


\subsection{Implementation details on the \gls{bosss} side}
\label{sec:impl-deta-glsb}

\gls{bosss} operators are usually part of a larger solver implementation, which
obviously makes it difficult to invoke them from \of{} in a straightforward
manner.
For that reason, a special \code{FixedOperators} class has been added to
\gls{bosss}. The purpose of this class is to provide wrapper methods for the
\gls{bosss} operators. These methods are in turn exposed to the \mono{} runtime, 
allowing them to be called from \of{}.

Currently, only a Laplace operator and the Cahn--Hilliard operator discussed
here have been implemented in this way, but adding additional operators from
\gls{bosss}'s vast collection of solvers is a straightforward task.

\subsection{Coupling of FVM and DG}
\label{sec:coupling-fvm-dg}

Naturally, both the \gls{fvm} and the \gls{dg} discretization operate on a
numerical mesh. However, whereas in the \gls{fvm} setting, the value of a
variable is constant across one cell, in \gls{dg}, the value of a variable is
described using an arbitrarily high-order polynomial.
Though this is somewhat of an oversimplification, \gls{fvm} can roughly be
viewed as a special case of \gls{dg} where the polynomials all have degree zero.
This is precisely the idea that allows translating between the two approaches
as shown in \cref{fig:flow}.

In the future, one may envision more sophisticated methods for coupling
\gls{fvm} and \gls{dg} being developed and implemented. Since \gls{fvm} requires
a finer numerical mesh than \gls{dg}, using the same mesh in both numerical
methods is not optimal. In the Cahn--Hilliard case, the resulting
performance penalty may be counteracted
by restricting the \gls{dg} calculation to a narrow band of cells close to the
interface, but in general, the requirement to use the same numerical mesh should
be eliminated. Obviously, this makes necessary the development of more
complicated synchronisation routines to communicate between the \gls{fvm} and
\gls{dg} worlds.



\section{Results}
\label{sec:results}

\subsection{One-dimensional profile}
\label{sec:one-dimens-prof}

\begin{figure}
  \centering
  \includegraphics[width=0.9\textwidth]{ 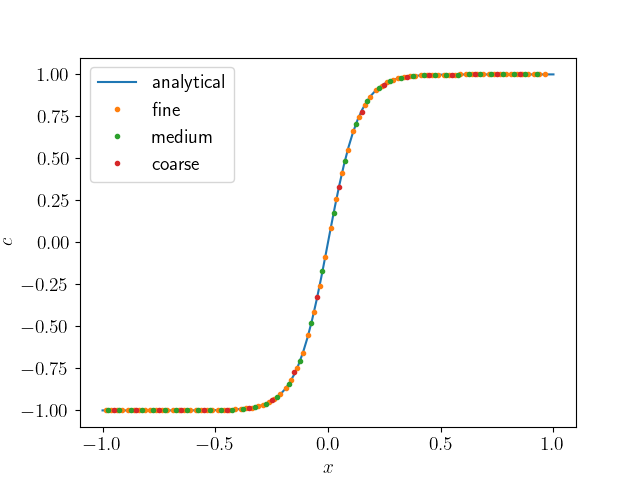 }
  \caption{One-dimensional profiles from three different mesh resolutions and
    the analytical solution.}
  \label{fig:oned}
\end{figure}

\begin{figure}
  \centering
  \includegraphics[width=0.9\textwidth]{ 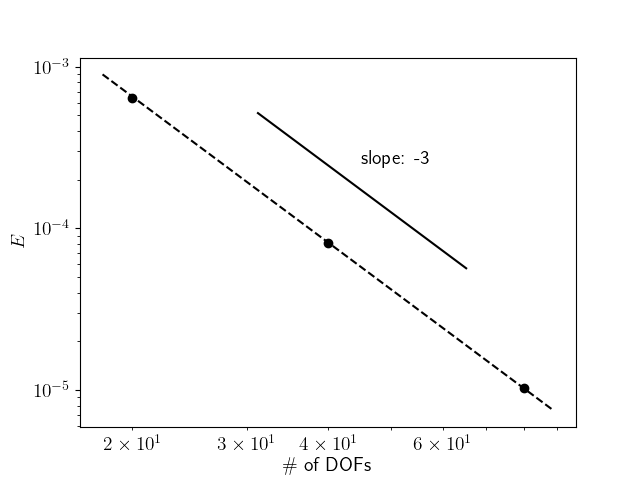 }
  \caption{Convergence plot of the calculations shown in \cref{fig:oned}.}
  \label{fig:onedconv}
\end{figure}

In order to verify the numerical implementation of the Cahn--Hilliard equation,
a simple one-dimensional test case is considered. On a domain \(x = [-1, 1]\),
\(c\) is initialized with
\begin{equation}
  \label{eq:3}
  c(x) = \sgn(x),
\end{equation}
which should converge to the analytical solution \parencite{cahn1958}
\begin{equation}
  \label{eq:4}
  c(x) = \tanh\left(\frac{x}{\sqrt{2 \gamma}}\right).
\end{equation}
The results of varying mesh sizes (20, 40 and 80 cells; all with DG-degree of
\(2\)) are
shown in \cref{fig:oned}.
Even the coarsest mesh shows no visible deviation from the analytical
solution, which indicates that we are well within the range of convergence.
As a confirmation, the errors are shown as a convergence plot in
\cref{fig:onedconv}.
A linear fit through the three points yields a convergence slope of \(2.99\),
which is in good agreement with the expected value of \(3\).

\subsection{Droplet in shear flow}
\label{sec:droplet-shear-flow}

\begin{figure}
  \centering
  \includegraphics[width=0.9\textwidth]{ 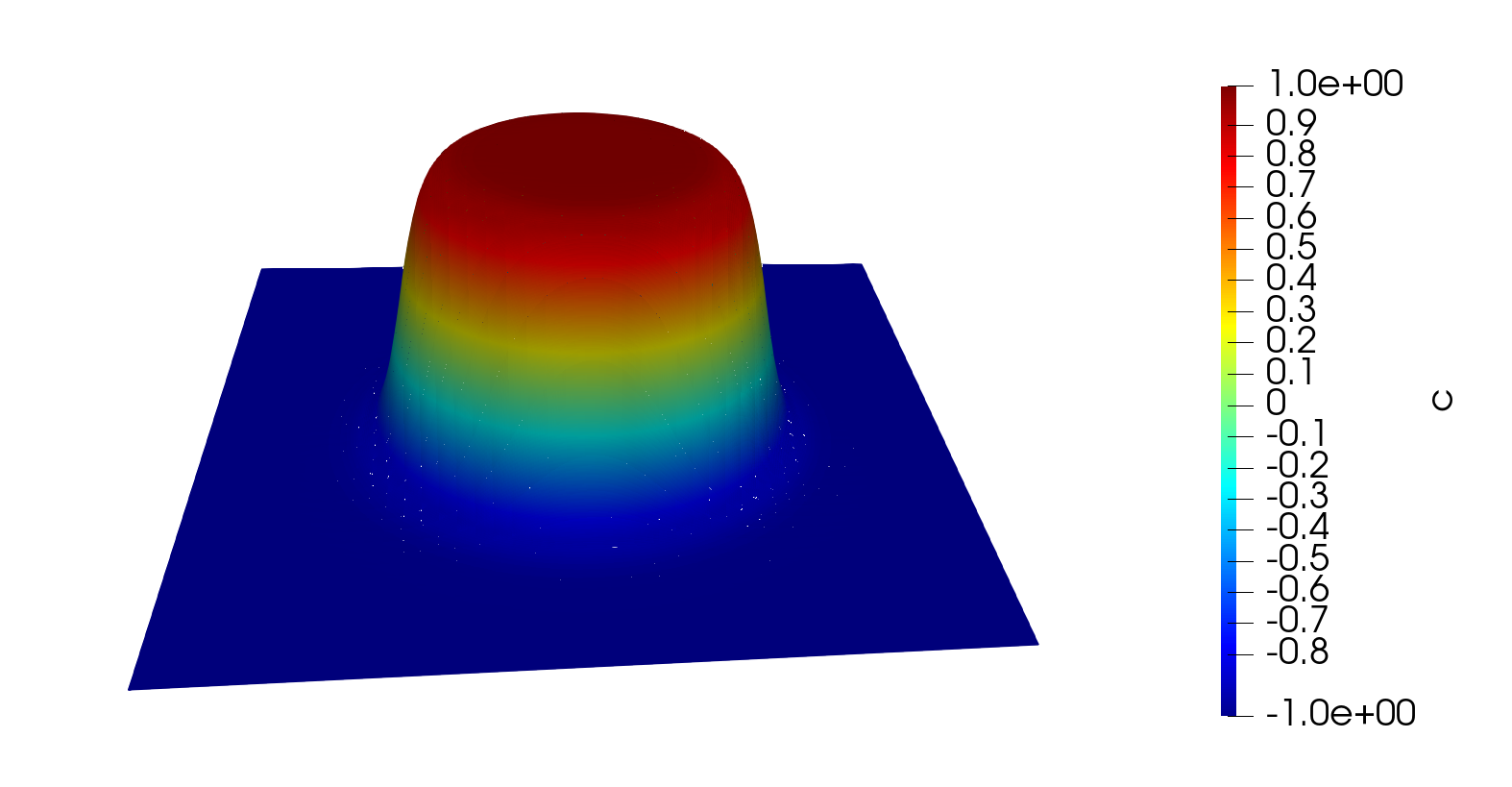 }
  \caption{Elevate plot of the order parameter field \(c\) at the final state of the droplet.}
  \label{fig:dropletElevate}
\end{figure}

As a two-dimensional test case for the Cahn--Hilliard equation and the coupling
to the momentum equation, we consider an initially round droplet that is
subjected to a shear flow.
The results shown in \cref{fig:dropletElevate}
exhibit a good qualitative behavior. Note again that any quantitative studies
will be the subject of future work.

\section{Conclusion and Outlook}
\label{sec:conclusion}

This article documents the usage of and some simple preliminary results obtained
using the \foamdg{} project, which provides the
possibility to conduct \gls{dg} calculations in \of{} by interfacing with
\gls{bosss}.
As showcased in the example discussed here, this interfacing can be done in a
highly modular manner, and allows for mixed approaches with some variables of
the problem being solved in \gls{fvm} by \of{} and others using \gls{dg} in
\gls{bosss}.
This method offers the attractive possibility of quickly prototyping the
enhancement of existing \of{} solvers with the increasingly popular high-order
\gls{dg} method.

In the future, more \gls{bosss} solvers could be made available to interface
with \of{}. While this process is straightforward, it does necessitate some
implementation work on the \gls{bosss} side.
Moreover, even though the current state works reasonably well, as mentioned in
\cref{sec:coupling-fvm-dg}, ideas for further efficiency and quality of life
improvements exist and could be implemented in the future.
Also, if the Cahn--Hilliard system is to be investigated further, some effort
has to be directed at understanding how the coupling between the equations is
affected by the combination of the different discretization schemes, in
particular with regard to the Eyre splitting of the potential term.

\ifdefined\review
\else
  \section*{Acknowledgements}

\noindent
The work of D. Klingenberg is funded by the Deutsche Forschungsgemeinschaft
(DFG, German Research Foundation) – Project-ID 265191195 – SFB 1194 and
supported by the 'Excellence Initiative' of the German Federal and State
Governments and the Graduate School of Computational Engineering at Technical
University Darmstadt.

\fi



\bibliographystyle{IEEEtran}

\bibliography{./latex_macros/bibliography.bib}

\end{document}